\documentclass[11pt]{amsart}
\usepackage{stmaryrd}
\usepackage{textcomp}
\usepackage{enumerate}
\usepackage{setspace}
\usepackage{amssymb}
\usepackage{amsthm}
\usepackage{amsmath}
\usepackage{graphicx}
\usepackage{extarrows}
\usepackage{marvosym}
\usepackage{empheq}
\usepackage{latexsym}
\usepackage{endnotes}
\usepackage{fontenc}
\usepackage{color}
\usepackage[usenames,dvipsnames]{xcolor}
\usepackage{comment}
\usepackage{url}
\usepackage{hyperref}
\usepackage{cleveref}
\usepackage{tikz}
\renewcommand{\epsilon}{\varepsilon}

\allowdisplaybreaks

\newtheorem{theorem}{Theorem}

\newtheorem{question}[theorem]{Question}
\newtheorem{lemma}[theorem]{Lemma}

\newtheorem*{claim*}{Claim}

\newtheorem{conj}[theorem]{Conjecture}

\newtheorem{observation}[theorem]{Observation}

\newtheorem{prop}[theorem]{Proposition}

\theoremstyle{definition}
\newtheorem{definition}[theorem]{Definition}

\theoremstyle{remark}
\newtheorem{remark}[theorem]{Remark}

\usepackage{thmtools}

\def\eps{\varepsilon}
\def\gam{\gamma}
\def\R{\mathbb{R}}
\def\sig{\sigma}

\renewcommand{\P}{\mathbb{P}}

\def\E{\mathbb{E}}

\def\Pt{\widetilde{P}}
\def\gt{\widetilde{\gamma}}
\def\tt{\tilde{t}}
\newcommand{\cE}{\mathcal{E}}

\newcommand{\OFG}{\operatorname{OFG}}
\newcommand{\tmix}{\tau_{\mathrm{mix}}}
\newcommand{\pit}{\widetilde{\pi}}
\newcommand{\TV}{\mathrm{TV}}

\newcommand{\sigsp}{\sig_{\mbox{sp}}}

\begin{document}

\title{On random locally flat-foldable origami}

\author{Thomas C. Hull} 
\address{Mathematics Department, Franklin and Marshall College}
\email{thomas.hull@fandm.edu}

\author{Marcus Michelen}
\address{Department of Mathematics, Statistics and Computer Science, University of Illinois Chicago.  Chicago, IL 60607}
\email{michelen.math@gmail.com}

\author{Corrine Yap}
\address{School of Mathematics, Georgia Institute of Technology, Atlanta, GA 30332}
\email{math@corrineyap.com}

\begin{abstract}
    We develop a theory of random flat-foldable origami. Given a crease pattern, we consider a uniformly random assignment of mountain and valley creases, conditioned on the assignment being flat-foldable at each vertex. A natural method to approximately sample from this distribution is via the face-flip Markov chain where one selects a face of the crease pattern uniformly at random and, if possible, flips all edges of that face from mountain to valley and vice-versa. We prove that this chain mixes rapidly for several natural families of origami tessellations---the square twist, the square grid, and the Miura-ori---as well as for the single-vertex crease pattern. We also compare local to global flat-foldability and show that on the square grid, a random locally flat-foldable configuration is exponentially unlikely to be globally flat-foldable.
\end{abstract}

\maketitle

\section{Introduction}\label{sec:introduction}

The folding of flat sheets of material into three-dimensional, sometimes highly intricate objects, commonly known as \textit{origami}, has seen a surge of interest for applications in physics \cite{Assis,Liu}, engineering \cite{Active,YouSci14}, architecture \cite{Maleczek2018}, and computer science \cite{GFALOP}. Origami can be modeled by a function $f:P\to\mathbb{R}^3$ on a compact domain $P\subset \R^2$ that is a continuous piecewise isometry; of particular interest are origami that \textit{fold flat}, meaning that the image $f(P)$ lies in the plane $\mathbb{R}^2$. Being able to fold flat allows for compact storage of the folded design, such as for transport and deploying of large solar panels into outer space \cite{YOrigami}.  

In this paper we develop a theory of \textit{random origami}, with the hope of answering questions such as ``With what probability will a randomly chosen crease pattern configuration be flat-foldable?" Our model is that of a uniformly random assignment of mountain/valley creases to a given crease pattern, conditioned on the crease assignment being locally flat-foldable.   

Already a natural question emerges: how does one study or (approximately) sample a uniformly random locally flat-foldable configuration?  Tools from the origami literature---Kawasaki's Theorem and Maekawa's Theorem detailed in \Cref{sec:background}---provide necessary and sufficient conditions for a configuration to be locally flat-foldable.  This allows one to recast the set of locally flat-foldable crease patterns under the umbrella of \emph{spin systems} from statistical physics and computer science.  Heuristically, this class of models is based on assigning configurations of some kind---in our case ``mountain'' and ``valley'' assignments to edges---under some family of local constraints.  In the cases of certain crease patterns, this shows an explicit connection between the topic of local flat-foldability and more well-studied questions such as proper coloring (see, for instance, \cite{Ginepro}).  

A broad class of problems studied across probability, statistical physics, and computer science is to understand how to approximately sample from a large system under local constraints.  Here, one typical approach is via \emph{Markov chain Monte Carlo}.  In our case, this amounts to the following: given one locally flat-foldable assignment, perform some small random change that transforms the assignment into another locally flat-foldable assignment.  It is often possible to choose these transition probabilities in a given way so that if one performs these random transitions over and over again, the resulting random configuration is approximately uniform over all valid configurations.  

In the setting of origami, a natural and easy-to-implement Markov chain is the \emph{face-flip} Markov chain: one takes a face of the crease pattern uniformly at random, checks if it is possible to flip all creases surrounding this face from mountain to valley and vice-versa, and performs this flip with probability $1/2$ if possible.  As a preliminary result, using basics from the theory of Markov chain mixing, we prove that this face-flip Markov chain converges to the uniform distribution on locally flat-foldable configurations for a wide family of crease patterns.  
Our primary settings are several lattices that are widely used throughout the origami literature, as well as single-vertex crease patterns. 
A main concern is to understand how long it takes before the Markov chain is approximately uniform. For this, we obtain bounds on the \emph{mixing time}; our main results are encompassed in the following theorem. 

\begin{theorem}\label{thm:main-lff}
    The face-flip Markov chain for locally flat-foldable mountain-valley assignments of the square twist, the square grid, the Miura-ori, and the single-vertex crease pattern all mix in polynomial time with respect to the number of faces.  
\end{theorem}
We prove each of the four cases of Theorem~\ref{thm:main-lff} separately, as Theorems~\ref{thm:mixing-square-twist}, \ref{thm:mixing-square-grid}, \ref{thm:mixing-miura}, \ref{prop:single-vertex-MC} respectively.

Polynomial-time mixing of the face-flip Markov chain immediately gives an approximate sampling algorithm that runs in polynomial time, which is called an {\em efficient} approximate sampler. We discuss this in more detail in \Cref{sec:background}.
For the single vertex and the square grid we additionally provide efficient exact samplers.
We also consider local versus global flat-foldability and prove that global flat-foldability is exponentially unlikely on the square grid.

\begin{theorem}\label{thm:globalff-square}
    Let $\sigma$ be a uniformly random locally flat-foldable configuration on the $m \times n$ square grid crease pattern. The probability that $\sigma$ is globally flat-foldable is $\exp[-cmn(1+o(1))]$ for some $c >0$. 
\end{theorem}

\subsection{Overview of Techniques}
Our primary technique for analyzing the face-flip Markov chain is to compare it to other Markov chains known to have polynomial mixing time. In the case of the first two lattice crease patterns---the square twist and the square grid---the comparison is quite straightforward. We use the perspective of face flips being equivalent to a lazy random walk on the origami flip graph (in other words, the Markov chain state diagram) and observe that the flip graph itself is isomorphic to a $d$-dimensional hypercube for some $d$. For the Miura-ori crease pattern, the route to comparison is more circuitous. Work by Ginepro and Hull \cite{Ginepro} showed a relationship between mountain-valley assignments on the Miura-ori crease pattern and proper vertex-colorings of grid graphs. Even further, a face flip on the Miura-ori crease pattern corresponds to a single update of \emph{Glauber dynamics}, a standard Markov chain in statistical physics that applies for proper colorings.  With this connection, we use powerful results about Glauber dynamics to obtain polynomial mixing for the face-flip Markov chain.

For single-vertex crease patterns, a  counting argument allows us to describe an exact sampler. To analyze the face-flip Markov chain, however, 
we pass to an induced chain on the configurations which have more mountain than valley creases. The induced chain can be compared to a type of fast-mixing card-shuffling Markov chain. 

Lastly, we consider local versus global flat-foldability on lattices.  It is simple to check whether a configuration is locally flat-foldable, but checking whether a configuration is globally flat-foldable is more complicated (see Section~\ref{sec:global} for more discussion on computational hardness results in this direction).  We show that on the square lattice, the probability that a random locally flat-foldable configuration is globally flat-foldable is exponentially small with respect to the volume of the lattice.

\subsection{Related Work}

Questions concerning local flat-foldability can be cast in the language of constraint satisfaction problems and statistical physics, specifically in the framework of \emph{spin systems}: each edge is assigned ``mountain'' or ``valley'' and there are local rules at each vertex to determine compatibility.  In particular, if one wishes to describe a uniformly random locally flat-foldable configuration as a spin system, one assigns energy $+\infty$ to each configuration that is not locally flat-foldable and energy $0$ to each configuration that is locally flat-foldable.  Via Kawasaki's Theorem, this can be checked by looking purely at each vertex individually. The perspective of interpreting locally flat-foldable configurations as a spin system is not new and indeed work by Assis \cite{Assis} draws connections between more classical statistical physics models and local flat-foldability on certain lattices (some of which we will also discuss in the present study).
Ginepro and Hull \cite{Ginepro} proved that the foldings of the \textit{Miura-ori} crease pattern (to be defined in Section~\ref{subsec:Miura}) are equivalent to the square-ice model. 
Nakajima used a spin model on random graphs to model the combinatorial problem of ordering the different layers of paper in flat-folded crease patterns \cite{nakajima1}, and developed more work in this vein for crease patterns that contain only one interior vertex \cite{nakajima2}. 

Our main focus is on properties of random locally flat-foldable configurations and in particular on sampling via Markov chains.  The use of Markov chains to sample spin systems is a classic topic at the intersection of statistical physics, probability theory, and computer science.  In particular, Markov chain Monte Carlo was introduced in a work by Metropolis--Rosenbluth--Rosenbluth--Teller--Teller \cite{metropolis1953equation} which sought to sample from a continuous-space variant of a spin system called the hard-disk model.  Analysis of Markov chains for spin systems remains a vibrant field, and we refer the reader to the following texts \cite{grimmett2018probability, levin2017markov} and some recent breakthroughs \cite{anari2021spectral,chen2021optimal} along with the references therein for more context.
There has also been quite a bit of research on the algorithmic complexity of origami over the past several decades, including several monographs \cite{GFALOP,Origametry,Ida,Uehara1}. We make references to related computational origami results throughout.

\subsection{Organization} The results of this paper lie at the intersection of origami, probability, and algorithms. Thus, we include background material in Section~\ref{sec:background} on each of the three areas; readers may let their experience dictate how much they choose to read. 
Section~\ref{subsec:origami} has preliminaries on flat origami and Section~\ref{subsec:prob} has preliminaries on Markov chains. Sections~\ref{subsec:gibbs} and \ref{subsec:face-flip} then contain preliminaries specific to uniformly random mountain valley assignments and the face-flip Markov chain. In Section~\ref{sec:lattice}, we prove \Cref{thm:main-lff} for the three types of aforementioned lattice crease patterns, and in Section~\ref{sec:single-vertex}, we prove the remaining case of the single-vertex crease patterns. In Section~\ref{sec:global}, we prove \Cref{thm:globalff-square}. In Section~\ref{sec:conclusion}, we discuss open problems and future directions.

\section{Preliminaries}\label{sec:background}

\subsection{Flat Origami}\label{subsec:origami}

Origami that can fold flat possess rich geometric and combinatorial structures. In what follows we adopt conventions of modeling flat origami presented in \cite{GFALOP,Origametry}. While we will not explicitly use many of the notions and notations defined below, they provide a foundation for posing questions about origami in a mathematically rigorous way. 

Given an origami $f:P\to\mathbb{R}^2$ that folds flat, we let the \textit{crease pattern} $C$ of $f$ be the subset of $P$ on which $f$ is non-differentiable.  It has been proven that $C$ forms a straight-line embedding of a graph on $P$ whose vertices always have even degree \cite{Robertson}. We call the edges of $C$ the \textit{creases} and the connected regions of $P$ between creases the \textit{faces} of $C$.  For any face $F$ of $C$, its image under the transformation $f$ either involves only translations and rotations and is thus isotropic or involves a reflection and is anisotropic. Faces of $C$ that share a common boundary crease have different isotropic/anisotropic orientations, and thus they define a proper 2-coloring $u_f:\{$faces of $C\}\to\{0,1\}$ of the faces of $C$. 

The folding map $f$ of an origami misses a crucial aspect of flat-folded paper: we need the paper $P$ to not \textit{self-intersect} under $f$ (i.e. we do not want the folded paper to be forced to penetrate itself in order to fold flat).  To capture this in our model, we need to describe the \textit{layer ordering} of $f$, denoted $\lambda_f:D\to\{\pm 1\}$, where $D$ is the subset of points $(p,q)\in P\times P$ such that both $p$ and $q$ are contained in the interior of the faces of $C$ and $f(p)=f(q)$. We say that $\lambda_f(p,q)=1$ (resp. $-1$) means $p$ is below $q$ (resp. $p$ is above $q$).  In order for $\lambda_f$ to guarantee that the flat origami map $f$ not self-intersect, it needs to satisfy certain \textit{non-crossing conditions} outlined in \cite{boxpleat,Origametry} but which do not need to be specified here. 

We may now define the \textit{mountain-valley (MV) assignment} $\sig:E(C)\to\{\pm 1\}$ of $f$ with layer ordering $\lambda_f$, where $E(C)$ is the set of edges of $C$, as follows: For any crease $e\in E(C)$ bordering faces $F_1$ and $F_2$ of $C$ and points $p\in F_1$ and $q\in F_2$ arbitrarily close but on opposite sides of $e$ with $f(p)=f(q)$, we have
$$\sigma(e)=1 \iff (u_f(F_1)=0\mbox{ and }\lambda_f(p,q)=1)\mbox{ or }(u_f(F_1)=1\mbox{ and }\lambda_f(p,q)=-1)\mbox{, and}$$
$$\sigma(e)=-1 \iff (u_f(F_1)=0\mbox{ and }\lambda_f(p,q)=-1)\mbox{ or }(u_f(F_1)=1\mbox{ and }\lambda_f(p,q)=1).$$
The above simply says that, relative to the orientation of the faces under the mapping $f$, $\sig(e)=1$ means the crease forms a \textit{valley} crease (and thus makes a $\vee$ shape when unfolded) and $\sig(e)=-1$ means that $e$ forms a \textit{mountain} crease (making a $\wedge$ shape when unfolded).  Layer orderings and MV assignments are strongly related; a layer ordering will generate an MV assignment on the creases, and an MV assignment on the creases will generate a layer ordering.  An MV assignment $\sig:E(C)\to\{\pm 1\}$ placed on a flat origami crease pattern $C$ is called \textit{valid} if it produces a layer ordering that avoids self-intersections of the paper.  

We often pay attention to how an origami crease pattern folds at each vertex individually, or \textit{locally}.  That is, given a crease pattern $C$ and an MV assignment $\sig$, we draw a sufficiently small disc $D_v$ around each vertex $v$ in a crease pattern $C$ so that the only vertex of $C$ in $D_v$ is $v$. Considering $C\cap D_v$ for each $v$ as separate crease patterns, we say that $C$ is \textit{locally flat-foldable} under $\sig$ (or alternatively, that $\sig$ is \textit{locally valid} for $C$) if each of the sub-crease patterns $C\cap D_v$ is flat-foldable under $\sig$. 

The following two fundamental results are very helpful in identifying and studying locally flat-foldable crease patterns (see \cite{GFALOP,Tom1,Origametry} for proofs):

\begin{theorem}[Kawasaki's Theorem]\label{thm:kawasaki}
    A crease pattern $(C,P)$ with only one vertex $v$ in the interior of $P$ and consecutive angles $\alpha_1,\ldots, \alpha_{2n}$ between the creases adjacent to $v$ is flat-foldable (meaning there exists a valid MV assignment $\sig$ for $C$) if and only if $\sum_{i=1}^{2n}(-1)^n\alpha_i=0$.
\end{theorem}

\begin{theorem}[Maekawa's Theorem]\label{thm:maekawa}
    A vertex $v$ in a flat-foldable crease pattern with valid MV assignment $\sig$ and creases $e_1,\ldots, e_{2n}$ adjacent to $v$ must satisfy
$\sum_{i=1}^{2n}\sig(e_i)=\pm 2$.
\end{theorem}

In general if a crease pattern has a valid MV assignment then we say that $C$ is \textit{globally flat-foldable}.  While local flat-foldability is easy to check via Kawasaki's Theorem, determining if a configuration is globally flat-foldabile is an NP-complete problem, even in the case where all the vertices of $C$ lie on the lattice $\mathbb{Z}^2$ and the angles between creases at the vertices are all multiples of $45^\circ$ (called \textit{box-pleated crease patterns}) \cite{boxpleat}.

Not every MV assignment is valid; in Section \ref{sec:global} we will see a simple example of a locally flat-foldable configuration that is not globally flat-foldable and use it to show that globally flat-foldable configurations are exponentially rare among locally flat-foldable configurations. 

\begin{figure}
    \centering
    \includegraphics[scale=.5]{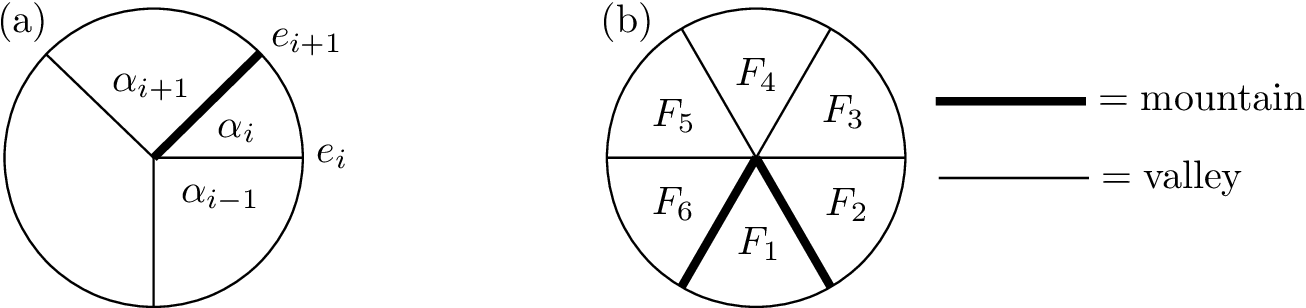}
    \caption{(a) A single-vertex crease pattern exhibiting the Big-Little-Big Theorem. (b) A degree-6 vertex with all equal angles between consecutive creases.}
    \label{fig:MVRulz}
\end{figure}

Consider the single-vertex crease pattern shown in Figure~\ref{fig:MVRulz}(a).  Creases $e_i$ and $e_{i+1}$ border the face with sector angle $\alpha_i$, which is strictly smaller that the sector angles $\alpha_{i-1}$ and $\alpha_{i+1}$ of its neighboring faces. In this situation, it is impossible for a valid MV assignment $\sig$ to have $\sig(e_i)=\sig(e_{i+1})$ because doing so in a flat folding of these creases would make the faces with bigger angles $\alpha_{i-1}$ and $\alpha_{i+1}$ cover the smaller $\alpha_i$ face on the same side of the paper, forcing a self-intersection. More formally, we have the following statement.

\begin{theorem}[Big-Little-Big Theorem]\label{thm:blb}
    If $\alpha_{i-1}, \alpha_i,\alpha_{i+1}$ are sector angles of consecutive faces around a flat-foldable vertex with $\alpha_i<\alpha_{i-1}, \alpha_{i+1}$, then a valid MV assignment $\sig$ must have $\sig(e_i)\ne\sig(e_{i+1})$, where $e_i$ and $e_{i+1}$ are the edges comprising angle $\alpha_i$.
\end{theorem}

\subsubsection{Face Flips}

We wish to study how a locally valid MV assignment evolves under random local changes; one of the most natural such operations is a face flip.

\begin{definition}Given a locally valid MV assignment $\sig$ and a face $F$ of our crease pattern $C$, we define the \emph{face flip of $F$} to be a new MV assignment $T_F(\sig)$ which equals $\sig$ for all creases that do not border $F$ and $T_F(\sig(e))=-\sig(e)$ for all creases $e$ that border $F$.  
\end{definition}

If $T_F(\sig)$ is also a locally valid MV assignment for $C$, then we say that face $F$ is \textit{flippable}, whereas a \textit{non-flippable face} of a locally valid MV assignment $\sig$ is a face $F$ for which $T_F(\sig)$ is not locally valid.

For example, in the locally valid MV assignment of the crease pattern shown in Figure~\ref{fig:MVRulz}(b) for the degree-6 vertex with all sector angles equal to $60^\circ$, the face $F_1$ is not flippable because doing so would make every crease a valley crease, violating Maekawa's Theorem.  But all of the other faces in this MV assignment are flippable.  In Figure~\ref{fig:MVRulz}(a), the face with sector angle $\alpha_i$ is flippable, but its neighbors with angles $\alpha_{i-1}$ and $\alpha_{i+1}$ are not flippable; in fact they are \textit{never flippable} faces, since for any locally valid MV assignment flipping the neighboring faces of the $\alpha_i$ face will violate the Big-Little-Big Theorem for creases $e_i$ and $e_{i+1}$.

Finally, for a given flat-foldable crease pattern $C$ we represent the state space of all locally valid MV assignments of $C$ as a graph.

\begin{definition}
The \emph{origami flip graph} of $C$, denoted $\OFG(C)$, is the graph whose vertices are all locally valid MV assignments of $C$ and where two vertices $\sig_1$ and $\sig_2$ form an edge if and only if there exists a face $F$ of $C$ such that $T_F(\sig_1)=\sig_2$.
\end{definition}

\subsection{Markov Chains}\label{subsec:prob}
Throughout we will make use of Markov chains on finite state spaces; we recall some classical notions from the area in this section and refer readers to~\cite{levin2017markov} for more background and context.  

\subsubsection{Basic definitions}

 A \emph{Markov chain} with \emph{state space} $\Omega$ and \emph{transition matrix} $P : \Omega \times \Omega \to \R$ is a sequence of random variables $X_0,X_1,\ldots$ so that for all $t \geq 0$ and values $x_j \in \Omega$ we have \begin{equation}\label{eq:markov-def}
    \P(X_{t + 1} = x_{t+1} \,|\, X_j = x_j \text{ for all }j \in \{0,1,\ldots,t\}) = \P(X_{t + 1} = x_{t+1} \,|\, X_t = x_t ) = P(x_t,x_{t+1})\,.
\end{equation}

When $P$ is viewed as an $|\Omega| \times |\Omega|$ matrix, \eqref{eq:markov-def} implies that the columns of $P$ sum to $1$.  The distribution of the trajectory $X_0,X_1,\ldots$ is determined uniquely by the distribution of $X_0$ along with the transition matrix $P$.  

A \emph{probability measure} on $\Omega$ is a function $\mu: \Omega \to \R$ with $\mu(x) \geq 0$ for all $x \in \Omega$ and $\sum_{x} \mu(x) = 1$.  For $A \subset \Omega$ we write $\mu(A) = \sum_{x \in A}\mu(x)$.  When $P$ is viewed as a matrix, it is convenient to view probability measures as row vectors and functions $f:\Omega \to \R$ as column vectors.  A {\em stationary distribution} for $P$ is a probability measure typically denoted $\pi$ such that $\pi P = \pi$, meaning for all $x \in \Omega$ we have $$ \pi(x) = \sum_{y} \pi(y) P(y,x)\,.$$ 

A sufficient condition for $\pi$ to be a stationary distribution of $P$ is if $\pi$ satisfies the {\em detailed balance equations},
\begin{equation} \label{eq:detailed-balance}
    \pi(x)P(x,y) = \pi(y)P(y,x)
\end{equation}  
for all $x, y \in \Omega$. In this case, we say $P$ is {\em reversible} with respect to $\pi$. 

A Markov chain is \emph{irreducible} if for each pair $x,y \in \Omega$ there is some $t \geq 1$ so that $P^t(x,y) > 0$, i.e.\ if it is always possible to get from any one state to any other by evolving according to the transition matrix $P$.
A Markov chain is {\em aperiodic} if $\gcd\{t \geq 1 : P^t(x,y) > 0\} = 1$ for all $x, y \in \Omega$. 
It will often be useful to assume that the Markov chain is \emph{lazy}, i.e. that for each $x$ we have $P(x,x) \geq 1/2$; laziness ensures that the Markov chain is aperiodic.

A fundamental theorem of Markov chains is that a Markov chain which is irreducible, aperiodic, and finite has a unique stationary distribution (see, e.g., \cite[Cor. 1.17]{levin2017markov}). 

\subsubsection{Mixing times}

Given two probability measures $\mu$ and $\nu$, define the \emph{total variation distance} via \begin{equation}
    \|\mu - \nu\|_{\TV} =  \max_{A \subset \Omega} |\mu(A) - \nu(A)|\,.
\end{equation}

If $P$ is irreducible and aperiodic with stationary distribution $\pi$, then we have the Markov convergence theorem, meaning that $\|P^t(x,\cdot) - \pi\|_{\mathrm{TV}}$ converges to zero (in fact, it does so exponentially quickly in $t$, see, e.g., \cite[Thm. 4.9]{levin2017markov}) . 

For algorithmic problems, one of the main interests is in the rate of convergence.  This is the motivation for the definition of the \emph{mixing time} of a Markov chain: 

\begin{definition}
We define the \emph{$\epsilon$-mixing time} to be \begin{equation}
    \tmix(\epsilon):= \min\{t : \max_{\mu} \| \mu P^t - \pi \|_{\TV} \leq \epsilon \} 
\end{equation}
where the maximum is over all probability measures $\mu$ on $\Omega$, and the {\em mixing time} to be $\tmix:= \tmix(1/4)$.
\end{definition}
We remark that the choice of $\frac14$ is by convention, as any choice of $\epsilon > 0$ will only give a constant factor difference. Specifically, $\tmix(\epsilon) \leq \lceil\log (1/\epsilon)\rceil \tmix(1/4)$.

Up to a constant of $2$, the mixing time can be reinterpreted using the following basic fact: 
\begin{equation} \label{eq:two-starting points}
    \sup_{\mu} \|\mu P^t - \pi \|_{\TV} \leq \sup_{x,y} \| P^t(x,\cdot) - P^t(y,\cdot)\|_{\TV} \leq 2 \sup_{\mu} \|\mu P^t - \pi \|_{\TV}\,.
\end{equation}

In the remainder of this subsection, we introduce some standard but slightly technical notions for giving upper bounds on mixing times; these techniques will be used in Section~\ref{sec:single-vertex} for single-vertex crease patterns.

All eigenvalues of a lazy reversible Markov chain $P$ lie in $[0,1]$, and so we may order them as $1 = \lambda_1 \geq \lambda_2 \geq \ldots \lambda_{|\Omega|} \geq 0$.  The \emph{spectral gap} is defined as $\gamma:= 1 - \lambda_2$.

By diagonalizing the matrix $P$ one obtains a bound on the mixing time (see e.g.\ \cite[Thm. 12.4]{levin2017markov}):

\begin{theorem}\label{thm:mixing-via-spectral}
    Let $P$ be a lazy, irreducible, reversible Markov chain with stationary distribution $\pi$.  Then $$\left(\frac{1}{\gamma} - 1\right)\log\left(\frac{1}{2\eps}\right) \leq \tmix(\eps) \leq \frac{1}{\gamma} \log\left(\frac{1}{\eps \pi_{\min}} \right)$$
    where $\pi_{\min} = \min\{ \pi(x) : x \in \Omega, \pi(x) > 0\}\,.$
\end{theorem}

One method for bounding the mixing time of a given Markov chain is to compare it to another Markov chain on the same state space, ideally one whose spectral gap is easier to analyze. To use comparison, we require the Dirichlet forms of both chains. 

Given a function $f : \Omega \to \mathbb R$, the \emph{Dirichlet form} $\cE$ of a reversible Markov chain $P$ on $\Omega$ with stationary distribution $\pi$ can be defined via 
\begin{equation}\label{eq:dirichlet}
\cE(f) = \frac{1}{2} \sum_{x,y \in \Omega} (f(x) - f(y))^2 \pi(x) P(x,y)\,.
\end{equation}

\begin{lemma}\label{lem:comparison}
    Let $P, \Pt$ be reversible Markov chains on $\Omega$ with stationary distributions $\pi, \pit$ respectively. If $\widetilde{\cE}(f) \leq \alpha \cE(f)$ for all $f$, then
    $$\gt \leq \max_{\sig \in \Omega} \frac{\pi(\sig)}{\pit(\sig)} \alpha \gamma\,.$$
\end{lemma}

One particular type of Markov chain we will use for comparison is an induced chain.

\begin{definition}\label{def:induced-chain}
Given a chain $(X_t)_t$ on $\Omega$ with transition matrix $P$ and $A \subset \Omega$, define $\tau_0 = \inf\{t : X_t \in A\}$ and iteratively define $\tau_k = \inf\{t > \tau_{k-1} : X_t \in A\}.$  Then the trajectory $(X_{\tau_0}, X_{\tau_1},\ldots)$ is a Markov chain on $A$ with transition matrix $P_A$ which we call the {\em induced chain on $A$}.
\end{definition}

For a lazy, reversible Markov chain let 
\begin{equation}\label{eq:hitting-time}
t_H(\alpha) = \max_{x \in\Omega, A \subset \Omega : \pi(A) \geq \alpha} \E_x[ \tau_A]
\end{equation}

where $\tau_A$ is the (random) time it takes for a trajectory started at $x$ to hit the set $A$.  The following theorem of Peres--Sousi relates this quantity to the mixing time.
\begin{theorem}[\cite{peres2015mixing}]\label{thm:hitting-time}
    For each $\alpha \in (0,1/2)$ there is a positive constant $C_\alpha$ so that 
    $$C_\alpha^{-1} t_H(\alpha) \leq \tmix \leq C_\alpha t_H(\alpha)\,.$$
\end{theorem}

\subsection{Uniformly random locally flat-foldable MV assignments}\label{subsec:gibbs}

The probability distribution given by taking a uniformly random locally flat-foldable configuration can be interpreted as a \emph{Gibbs measure}.  This means that if $\sigma$ is a uniformly random locally flat-foldable configuration then there is some energy function $H: \{\pm 1\}^{E(C)} \to [0,+\infty]$ called the \emph{Hamiltonian} for which we have that $$\P(\sigma = \tau) = \frac{1}{Z} e^{-H(\tau)} \quad \text{ with } Z = \sum_{\tau} e^{-H(\tau)}\,.$$

In the case of local flat-foldability, the Hamiltonian $H$ will simply check that Maekawa's theorem (\ref{thm:maekawa}) holds at each vertex.  We now explicitly describe the uniformly random locally flat-foldable configuration in terms of a Gibbs measure.

For a given vertex $v \in V(C)$, let $N_v \subset E(C)$ denote the set of edges incident to $v$.  Define $\phi_v : \{\pm 1\}^{E(C)} \to \{0,1\}$ by $$\phi_v(\tau) = \begin{cases}
        1 &\text{ if }\tau \text{ is locally flat-foldable at }v \\
        0 &\text{ otherwise}
    \end{cases}$$
and note that we may in fact view $\phi_v$ as a function on $\{\pm 1\}^{N_v}$ since it depends only on the values of $\tau$ on $N_v$.
We note also that $\tau$ is locally flat-foldable if and only if $\prod_{v \in V(C)} \phi_v(\tau) = 1$.  
Defining $$Z = \sum_{\tau \in \{\pm 1\}^{E(C)}} f_v(\tau)$$ we see that $Z$ is the total number of locally flat-foldable configurations and that $$\P(\sigma = \tau) = \frac{1}{Z} \prod_{v \in V(C)} \phi_v(\tau|_{N_v})\,.$$ 

One important property of Gibbs measures is that they satisfy the following domain Markov property: taking a uniformly random configuration and resampling a subset, conditioned on its boundary, remains uniformly random.  In particular we have the following basic fact, which follows from \cite[Theorem 7.12]{grimmett2018probability}:

\begin{lemma}\label{lem:DLR}
Let $C$ be a crease pattern and $\sigma : E(C) \to \{\pm 1\}$ be a locally flat-foldable MV assignment chosen uniformly at random.  For a collection $S \subset E(C)$, let $N(S) \subset E(C)$ denote the set of edges in $E(C)$ that are incident to some edge in $S$.  Then for any set $T \subset E(C) \setminus S$ and configurations $\tau_S,\tau_{T}$ with $\P(\sigma|_T = \tau_{T}) > 0$ we have $$\P(\sigma|_S = \tau_S \,|\, \sigma|_{T} = \tau_T) = \P(\sigma|_S = \tau_S \,|\, \sigma|_{T \cap N(S)} = \tau_T|_{T \cap N(S)})\,.$$
\end{lemma}

\subsection{The Face-Flip Markov Chain}\label{subsec:face-flip}

We now introduce the main character of this paper, the {\em face-flip Markov chain}:
\begin{enumerate}
    \item Given an MV assignment of $C$, choose a face $F$ uniformly at random.
    \item If $F$ is non-flippable, do nothing. If $F$ is flippable, do nothing with probability $\frac12$ and flip $F$ with probability $\frac12$. 
\end{enumerate}
Our goal will be to analyze the mixing time of this Markov chain on specific crease patterns, but we begin by making some general observations about this random process. By the discussion in the previous subsection, this chain has a unique stationary distribution which we can easily check must be uniform.

\begin{lemma}
    Given a crease pattern $C$, the face-flip Markov chain on $C$ is a reversible Markov chain whose stationary distribution is uniform over all locally flat-foldable MV assignments of $C$.
\end{lemma}
\begin{proof}
    Letting $\Omega$ be the set of all locally flat-foldable MV assignments, a face flip is an involution so $P(\sig, \sig') = P(\sig', \sig)$ for all $\sig, \sig' \in \Omega$. Thus, the uniform distribution satisfies the detailed balance equations \eqref{eq:detailed-balance}.
\end{proof}

From this, we see that the face-flip Markov chain gives an approximate sampling algorithm for uniformly random locally flat-foldable MV assignments: simply run the face-flip Markov chain for $t$ steps where $t \geq \tmix$. By definition, the resulting distribution is close in total variation distance to the uniform distribution (hence an ``approximate'' and not an exact sampler).

Another way we can view the face-flip Markov chain on a crease pattern $C$ with $n$ faces is as a type of {\em lazy random walk} on $\OFG(C)$: starting from a vertex $v \in \OFG(C)$, for each neighbor $u$ of $v$, move to $u$ with probability $\frac{1}{2n}$; with the remaining probability, stay at $v$.

We can also consider random walks on general graphs, and this will prove useful in analyzing the square twist and square grid crease patterns. We will consider the particular graph of the {\em $d$-dimensional hypercube}, denoted $Q_d$. There are multiple equivalent descriptions of $Q_d$; the one most useful to us is where we think of the vertices as binary strings of length $d$ and edges consisting of pairs of vertices with Hamming distance 1. More precisely, $V(Q_d) = \{0, 1\}^d$ and $\{(u_1, u_2, \dots, u_d), (v_1, v_2, \dots, v_d)\} \in E(Q_d)$ if and only if there is a unique $i$ such that $u_i \neq v_i$. 

Another perspective we will use is that of a Cayley graph. Given a group $(G,\cdot)$ and generating set $S$, the {\em Cayley graph} $\Gamma = \Gamma(G, S)$ has vertex set $V(\Gamma) = G$ and $gh \in E(\Gamma)$ if and only if there exists $s \in S$ such that $gs = h$. We can naturally identify $Q_d$ with the Cayley graph of $\mathbb F_2^d$ with generating set $S = \{e_i\}_{i=1}^d$ where $e_i = \{0, \dots, 0, 1, 0, \dots, 0\}$ is the $i$th standard basis vector, having a 1 in the $i$th coordinate.

Observe that if $\OFG(C)$ is isomorphic to $Q_d$ for some $d$, then the face-flip Markov chain is a {\em simple} lazy random walk, meaning at every vertex $v$, there is a $\frac12$ probability of staying at $v$ and a $\frac{1}{2\deg(v)}$ probability of moving to each neighbor. Furthermore, as $\deg(v) = d$ for all $v \in V(Q_d)$, we can check via the detailed balance equations \eqref{eq:detailed-balance} that the unique stationary distribution is uniform over all $V$.

The mixing time of the lazy simple random walk on $Q_d$ is well-studied (see, e.g., \cite{levin2017markov}).

\begin{theorem}\label{thm:hypercube-mixing}
    For the lazy simple random walk on $Q_d$, 
    $$\tmix = O(d \log d)\,.$$
\end{theorem}

As a remark, this result follows from analysis of the coupon collector problem, which asks: If there are $n$ coupons and each week a person receives a coupon independently and uniformly at random, how many weeks does it take to collect all coupons?  The expected number of weeks is precisely  the $n$th harmonic number, asymptotically equal to $n \log n + O(1)$. The time it takes for the random walk on $Q_d$ to mix can be bounded by the time it takes to select each of the $d$ coordinates at least once---by coupon collector, this is $O(d \log d)$.

\section{Lattice crease patterns}\label{sec:lattice}

In this section, we discuss uniformly random locally flat-foldable crease patterns on various lattices and utilize known results about each origami tessellation in order to prove results about mixing times of the face-flip Markov chain, as promised in \Cref{thm:main-lff}.

\subsection{Square twist}\label{sec:sqtwist}

\begin{figure}
    \centering
    \includegraphics[scale=.55]{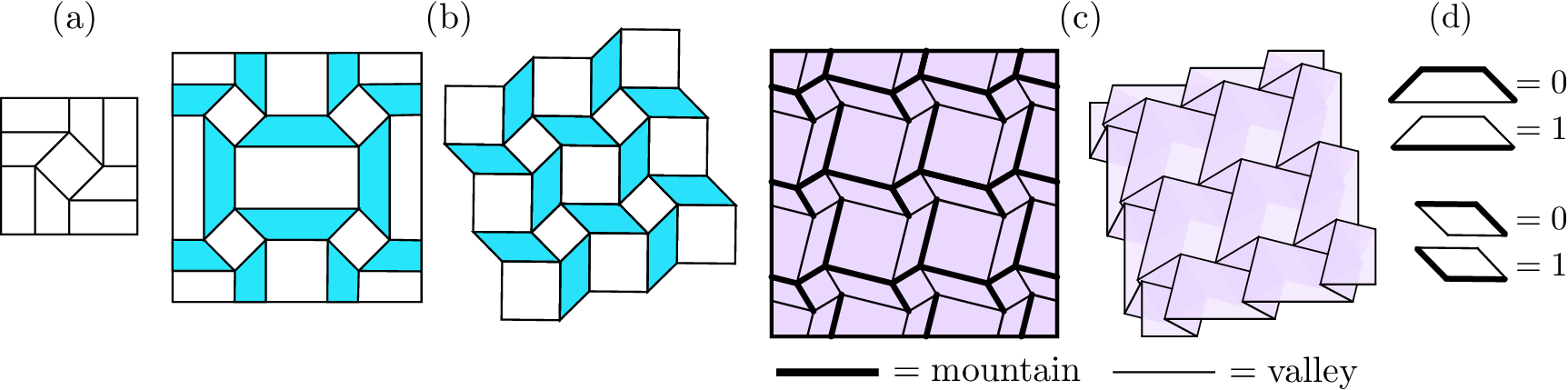}
    \caption{(a) The square twist crease pattern. (b) Two ways to tessellate square twists, with isotropic/anisotropic 2-colorings shown. (c) A valid MV assignment for a square twist tessellation, with an illustration of the folded image (made using \cite{TessLang}). (d) The possible MV assignments for trapezoid and parallelogram faces.}
    \label{fig:sqtwist}
\end{figure}

The \textit{square twist} crease pattern is shown in Figure~\ref{fig:sqtwist}(a).  It is used as a building block in many \textit{origami tessellations}, which are crease patterns that form a tiling of the plane if extended indefinitely \cite{Gjerde}. Square twists can be tiled in a grid either with alternating chirality or with all the same chirality, as shown in Figure~\ref{fig:sqtwist}(b). An example of a valid MV assignment for such a tessellation, along with a see-through simulation of its folded image (made using Lang's \textit{Tessellatica} software \cite{TessLang}) is in Figure~\ref{fig:sqtwist}(c). (The pattern on the right in Figure~\ref{fig:sqtwist}(b) is an example of Barreto's Mars, discussed in \cite{Assis}.) In addition to being artistically pleasing, square twists exhibit interesting mechanical properties when folded and unfolded, making them of interest for applications \cite{Silverberg}.

The proper 2-face-coloring obtained from square twist tessellations, which partitions the faces into rectangular and non-rectangular classes, is especially useful for describing their origami flip graphs, which then allow us to describe the probability space of all valid MV assignments. Figure~\ref{fig:sqtwist}(b) shows these 2-colorings.  

\begin{prop}
    The face-flip graph of a square twist tessellation $C$ is isomorphic to the $d$-dimensional hypercube where $d$ is the number of non-rectangular faces of $C$.
\end{prop}

\begin{proof}
Note that by the Big-Little-Big Theorem (\ref{thm:blb}), the square and rectangular faces of $C$ are never flippable (since flipping one from a valid MV assignment will make a $45^\circ$ angle have two Ms or two Vs). For the non-rectangular faces we observe the following.

\begin{itemize}
    \item They are either trapezoids or parallelograms (see Figure~\ref{fig:sqtwist}(b)).
    \item Their $45^\circ$ angles must be bordered by a mountain and a valley (by Big-Little-Big (\ref{thm:blb})).
    \item Their $135^\circ$ angles must be bordered by two mountains or two valleys (since each vertex has degree 4 and via Maekawa's Theorem (\ref{thm:maekawa})).
\end{itemize}

Therefore in a valid MV assignment of a square twist tessellation, the non-rectangular faces must have mountains and valleys in one of the arrangements shown in Figure~\ref{fig:sqtwist}(d). This implies that these faces can only have MV assignments that come in flippable pairs. Thus the non-rectangular faces
are always flippable (under all valid MV assignments), and flipping any non-rectangular face does not affect the MV assignment of any other non-rectangular face. 

This means that if $d$ is the number of non-rectangular faces in the crease pattern $C$, we can describe any valid MV assignment $\sigma_1$ of $C$ as a vector $v_1\in \{0,1\}^d$ where the $i$th coordinate of $v_1$ describes which of the two possible MV assignments (see Figure~\ref{fig:sqtwist}(d)) of the $i$th non-rectangular face we have in $\sigma_1$. If $\sigma_2$ is another valid MV assignment of $C$ and is described by vector $v_2$, then a sequence of faces to flip from $\sigma_1$ to $\sigma_2$ will be those whose coordinates in $v_1$ and $v_2$ are not equal. This proves that $\OFG(C)$ is connected and isomorphic to the hypercube $Q_d$.
\end{proof}

As discussed in the preliminaries, the mixing time of the simple lazy random walk on the $d$-dimensional hypercube is well-known---by \Cref{thm:hypercube-mixing}, we have the following. 

\begin{theorem}\label{thm:mixing-square-twist}
    The mixing time of the face-flip Markov chain of the square twist with $n$ faces is $O(n \log n)$.
\end{theorem}

\subsection{Square grid}

The \textit{square grid} crease pattern $S_{m,n}$ consists of $m$ rows and $n$ columns of square faces. They are the subject of the \textit{stamp folding problem}, a problem with some history \cite{koehler,lunnon1}, but such prior work often strives to count the different layer orderings possible for a given MV assignment.  In the present work, we focus only on distinct MV assignments.

We first list some observations about the square grid; see \cite{faceflips} for proofs.
\begin{itemize}
    \item Given a locally flat-foldable configuration, any face can be flipped and the result will still be locally flat-foldable. 
    \item Given two locally flat-foldable configurations, there exists a unique set of face flips that transforms one into the other.
    \item Face flips commute.
\end{itemize}

As with the square twist, we compare the square grid origami flip graph to a hypercube of some dimension.  One can interpret the hypercube as the Cayley graph of $\mathbb{F}_2^d$ with the usual generators $\{e_j\}$.  Conversely, we show that the flip graph is isomorphic to a Cayley graph of a quotient of $\mathbb{F}_2^d$ by the all ones vector.  This is the same as taking the hypercube and identifying each pair of antipodal points (say, $\vec{0}$ and $\vec{1}$). 

\begin{prop}
    The face flip graph of an $n$-face square grid is isomorphic to the Cayley graph $\Gamma$ of $\mathbb{F}_2^n / \langle 1,1,\ldots, 1\rangle$ generated by the set of standard basis vectors $\{e_j\}_{j=1}^n$.
\end{prop}

\begin{proof} We use our previous observations about the square grid and construct an isomorphism $\phi$ between the face-flip graph $F$ and $\Gamma$ as follows: fix a locally flat-foldable configuration $\sig$ and let $\phi(\sig) = \vec{0}$. Order the faces of the square grid $f_1, f_2, \dots, f_n$ and associate each element $\sig'$ of $\Omega$ with $\phi(\sig') = \vec{v} \in \{0, 1\}^n$ such that $v_i = 1$ if and only if the set of faces transforming $\sig$ to $\sig'$ contains $f_i$ (we know this set is defined and unique by our earlier observations). It is clear that $\phi^{-1}(\vec{1}) = \sig$ since flipping every single face will reverse each crease exactly twice. 

To see that $\phi$ is a bijection, consider $\vec{v} \in \{0,1\}^n$ and its corresponding set $S$ of faces which when flipped on $\sig$ result in a configuration $\sig'$. In order to achieve $\sig'$ by a different set of face flips $T$, we must have for every crease $e$ that the parity of the number of faces adjacent to $e$ is the same in $S$ as it is in $T$. If $S \neq T$, then this determines the entire configuration: $T$ must be $S^c$. 
\end{proof}

By the same coupon collector argument used to prove \Cref{thm:hypercube-mixing}, we conclude the following. 

\begin{theorem}\label{thm:mixing-square-grid}
    The mixing time of the face-flip Markov chain of the square grid with $n$ faces is $O(n \log n)$.
\end{theorem}

\subsection{Miura-ori}\label{subsec:Miura}

\begin{figure}
    \centering
    \includegraphics[width=\linewidth]{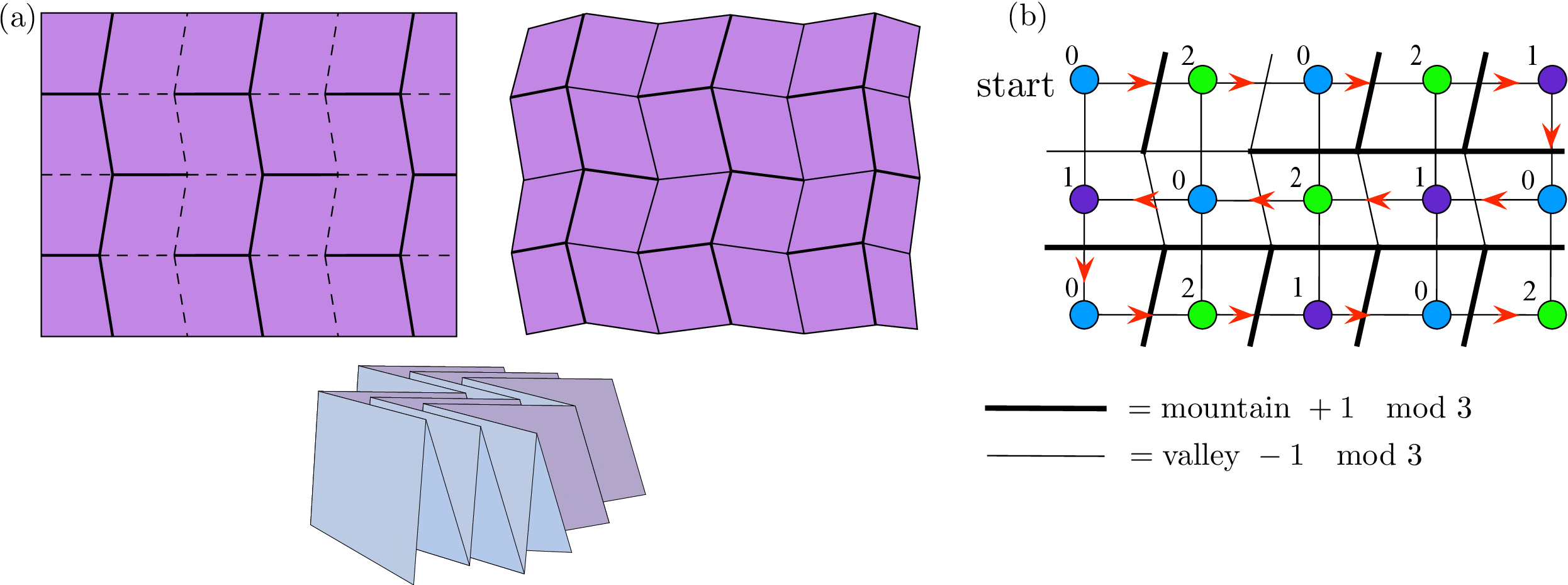}
    \caption{(a) A $4\times 6$ example of the Miura-ori crease pattern, with  MV assignment that folds up in the standard way. (b) Illustrating the bijection between locally-valid MV assignments of an $m\times n$ Miura-ori and proper 3-vertex colorings of the $m\times n$ grid graph (with one vertex pre-colored).}
    \label{fig:Miura}
\end{figure}

The \textit{Miura-ori} crease pattern is an origami tessellation whose faces are all congruent parallelograms, arranged as in Figure~\ref{fig:Miura}(a).  It is possibly the most well-studied origami crease pattern, as it was originally designed by Japanese astrophysicist Koryo Miura in the 1970s as a way to collapse large solar panel arrays and then transport and deploy them into outer space \cite{YOrigami}, and has since attracted attention for its interesting mechanics and use as a tunable metamaterial \cite{Silverberg1}.

To analyze the face-flip Markov chain, we pass to a well-studied graph coloring process. The following two lemmas give the relationship between Miura-ori crease patterns and vertex colorings.

\begin{lemma}[\cite{Ginepro}]\label{lem:miura-3-color}
    There exists a bijection $\psi$ between the set of locally flat-foldable configurations on the $m \times n$ Miura-ori crease pattern and the set of proper 3-vertex-colorings of the $m \times n$ grid graph with one vertex pre-colored.
\end{lemma}

\begin{lemma}[{\cite[Lemma 4.1]{faceflips}}]\label{lem:miura-local}
    Let $\sig, \sig'$ be locally valid Miura-ori configurations. Then $\sig' = T_F(\sig)$ if and only if $\psi(\sig)$ and $\psi(\sig')$ differ precisely at the vertex corresponding to $F$.
\end{lemma}

Thus, to bound the mixing time of face flips for Miura-ori crease patterns, it suffices to bound the mixing time for an analogous random process on proper $3$-vertex-colorings of grid graphs, where at each step a vertex is chosen uniformly at random and assigned a possible color (that preserves the properness of the 3-coloring) uniformly at random. This process is known as the {\em Glauber dynamics} for proper colorings and is a type of Markov chain that is often studied in statistical physics.

\begin{lemma}[{\cite[Theorem 3.1]{Goldberg}}]\label{lem:glauber-mix}
    The Glauber dynamics on proper 3-vertex-colorings of the $m \times n$ grid mix in time $O(m^4 n^9)$.
\end{lemma}

We then immediately obtain our main result.

\begin{theorem}\label{thm:mixing-miura}
    The mixing time of the face-flip Markov chain of the $m \times n$ Miura-ori crease pattern is $O(m^4n^9)$.
\end{theorem}

\section{Single-vertex crease patterns}\label{sec:single-vertex}

Flat-foldable crease patterns that consist of a single vertex in the paper's interior, which we denoted $C\cap D_v$ in Section~\ref{subsec:origami}, form the building blocks of flat-foldable crease patterns. 
Kawasaki's Theorem (\ref{thm:kawasaki}) provides us with an easy-to-check necessary and sufficient condition to determine whether a single-vertex crease pattern is flat-foldable, but it does not provide an MV assignment to fold the vertex flat. Further, Maekawa's Theorem (\ref{thm:maekawa}) is only a necessary condition. Hull \cite{Tom2} and Justin \cite{Justin} describe how to recursively enumerate all valid MV assignments of a flat-foldable single-vertex crease pattern, providing the tools needed to determine single-vertex flat-foldability in linear time \cite{GFALOP}.  Given all this, the structure of valid MV assignments for single-vertex crease patterns under face flips is surprisingly complex.  Origami flip graphs for single-vertex flat-foldable crease patterns are not well understood, even in the case where all the sector angles between adjacent creases are congruent \cite{JGAA-605}.

Along these lines, we let $C_{2n}$ be the crease pattern consisting of a single vertex with $2n$ equally spaced creases incident to it.  An example of $C_6$ is shown in Figure~\ref{fig:MVRulz}(b). Set $\Omega$ to be the collection of locally flat-foldable MV configurations of $C_{2n}.$

\subsection{Perfect sampler}

It is not too difficult to describe an exact sampling algorithm for this model that runs in linear time.

\begin{prop}\label{pr:single-vertex-uniform}
    There is a perfect sampler for the uniform distribution on $\Omega$ that runs in time $O(n).$   
\end{prop}

The idea will be to generate each valid MV assignment iteratively, conditioned on the previously revealed creases.  Label the creases $1,2,\ldots,2n$. Then a configuration $\sigma \in \Omega$ may be viewed as a function $\sigma:[2n] \to \{-1,1\}$ with $\sum_j \sigma(j) \in \{-2,2\}\,.$  Writing $\P$ for the uniform probability measure on $\Omega$, we begin by computing the marginal probabilities.

\begin{lemma}\label{fact:single-marginal}
    For each $j \in [2n]$ and $r \in [-2n, 2n]$, we have $$\P( \sigma(j) = +1 \,|\, \sigma(1),\ldots,\sigma(j-1)) = \P(\sigma(j) = +1 \,|\, \sum_{i < j} \sigma(i) = r) = \frac{\binom{2n - j}{(2n - j-r+1)/2} + \binom{2n - j}{(2n - j -r-3)/2}}{\binom{2n-j+1}{(2n-j-r+3)/2}+\binom{2n-j+1}{(2n-j-r-1)/2}}.
    $$
\end{lemma}

\begin{proof} 
    First note that there are exactly $\binom{m}{(m + k)/2}$ elements of $\{-1,1\}^m$ that sum to $k$, since each such sequence must contain exactly $k+(m-k)/2 = (m + k)/2$ many $1$'s and $(m - k)/2$ many $-1$'s. Given $\sigma(1), \dots, \sigma(j-1)$, by Maekawa's Theorem (\ref{thm:maekawa}) the valid MV assignments of $\{\sig(j), \dots, \sig(2n)\}$ are in bijection with sequences of length $2n-(j-1)$ that sum to either $2-r$ or $-2-r$ where $r = \sum_{i<j} \sig(i)$; this implies the first equality in the claim. 
    For the second equality, we compute that the number of choices for $\{\sig(j), \dots, \sig(2n)\}$ is thus 
    $$\binom{2n-j+1}{(2n-j-r+3)/2}+\binom{2n-j+1}{(2n-j-r-1)/2}.$$
      
    The total number of such choices that additionally have $\sigma(j) = + 1$ is the number of sequences of length $2n - j$ that sum to either $1 - r$ or $-3 - r$, thus giving the count of 
    $$\binom{2n - j}{(2n - j-r+1)/2} + \binom{2n - j}{(2n - j -r-3)/2}.$$
    Taking the ratio of the previous two displayed equations completes the proof.
\end{proof}

\begin{proof}[Proof of Proposition \ref{pr:single-vertex-uniform}]
    The steps of the perfect sampler are as follows: for each $j \in [2n]$ compute the marginal probability for $\P(\sigma(j) = +1 \,|\, \sigma(1),\ldots,\sigma(j-1))$ using \Cref{fact:single-marginal}.  Choose $\sigma(j) = +1$ with that probability and choose $\sigma(j)=-1$ otherwise.  Perform this operation for $j = 1,2,\ldots,2n$.  Each step may be done in constant time and there are $O(n)$ steps total.
\end{proof}

\subsection{Markov chain}

The main result of this section is the following.

\begin{theorem}\label{prop:single-vertex-MC}
    The mixing time of the face-flip Markov chain on $C_{2n}$ is $O(n^{5}\log n).$
\end{theorem}

Before proving this, we introduce several key ingredients.
Let $P$ denote the transition matrix for the face-flip random walk. From Maekawa's Theorem (\ref{thm:maekawa}) and its converse for single-vertex crease patterns---shown in \cite{faceflips}---we have the following observation.
\begin{observation}\label{obs:number-mountains}
    $\Omega$ consists of all configurations with either $n+1$ or $n-1$ mountain assignments.
\end{observation}

To analyze $P$, we categorize the faces of a given MV assignment: say a face is MM if bordered by two mountain creases, VV if bordered by two valley creases, and mixed otherwise. Let $\Omega_{M} \subset \Omega$ be the set of configurations with $n+1$ mountain creases. Observe that for $\sigma \in \Omega_M$, flipping an MM face results in a configuration in $\Omega_M^c$ and similarly, flipping a VV face puts us back in $\Omega_M$. We wish to analyze $\Omega_M$ by itself (our arguments will analogously apply to $\Omega_M^c$). To that end, let $\pi_M$ be the uniform measure on $\Omega_M$. Recalling \Cref{def:induced-chain}, let $P_{M}$ be the transition matrix of the induced Markov chain where we restrict to the set $\Omega_{M}$.

To bound the mixing time of $P_M$, we compare this process to the well-studied {\em adjacent transposition shuffle} (see, e.g., \cite[16.1]{levin2017markov}). The latter is defined on permutations of $\{1, \dots, n\}$ and consists of choosing an adjacent pair $(i,i+1)$ uniformly at random (where $1 \leq i \leq n-1$) and swapping their values. We can view the face-flip process similarly: 
arbitrarily choose a crease to label 1 and label the remaining creases $2, \dots, 2n$ proceeding clockwise from 1. Let $\Pt_M$ be the face-flip chain on $\Omega_M$ but disallowing the flip $(1,2n)$.

Applying $\Pt_M$ to $\sigma \in \Omega_M$ is equivalent to assigning an arbitrary labeling $M_1, M_2, \dots, M_{n+1}$, $V_{n+2}, V_{n+3}, \dots, V_{2n}$ to the image of $\sigma$, running the adjacent transposition shuffle on this labeling, and then removing the labels from the output. Thus, we can see that if the adjacent transposition shuffle mixes rapidly, then $\Pt_M$ must also mix rapidly. We will show by comparison that this implies $P_M$ also mixes rapidly and from this deduce that $P$ mixes rapidly.

The mixing time of the adjacent transposition shuffle is well-studied.
\begin{lemma}[\cite{wilson2004mixing}]\label{lem:adj-transposition}
    The adjacent transposition shuffle on $[n]$ mixes in time $O(n^3 \log n)$.
\end{lemma}

With this, we proceed to proving our main result of this section.

\begin{proof}[Proof of \Cref{prop:single-vertex-MC}]
    Let $\gam, \gam_M, \gt_M$ be the spectral gaps and $t, t_M, \tt_M$ be the mixing times for $P, P_M, \Pt_M$ respectively.
    By \Cref{lem:adj-transposition}, we have $\tt_M = O(n^3 \log n)$ and so \Cref{thm:mixing-via-spectral} implies $\frac{1}{\gt_M} = O(n^3 \log n)$. Applying \Cref{lem:comparison} with $\alpha = 1$ gives $\gt_M \leq \gam_M$ and thus by \Cref{thm:mixing-via-spectral} we have that 
    $$t_M \leq \frac{1}{\gt_M} \log (2^{2n}) \leq O(n^4 \log n)$$
    
    To use this to bound $t$ we analyze the expected hitting time defined by \eqref{eq:hitting-time}. Let $A \subset \Omega$ with $\pi(A) \geq 1/4.$  We have at least one of $\pi(A \cap \Omega_M) \geq 1/8$ or $\pi(A \cap \Omega_V) \geq 1/8$; assume without loss of generality that $\pi(A \cap \Omega_M) \geq 1/8$ and write $A_M = A \cap \Omega_M$.  For any $x \in \Omega$, consider a trajectory $(x = X_0,X_1,\ldots)$ evolving according to $P$ and let $(X_{\tau_0}, X_{\tau_1}, \dots)$ be the chain induced on $\Omega_M$ with transitions $P_M$ (recall that we iteratively define $\tau_k = \inf\{t > \tau_{k-1} : X_t \in \Omega_M\}$). Let $\tau_{A_M}$ be the hitting time of $A_M$ in $(X_t)_{t \geq 0}$, and let $\tau_{A_M}^{M}$ be the hitting time of $A_M$ in the induced chain. Note that $\tau_A \leq \tau_{A_M}$, so in order to upper bound $\E_x[\tau_{A}]$, it suffices to bound $\E_x[\tau_{A_M}]$.

    Observe that for each $n \geq 1$ and for $C \geq 1$ constant, we can write
    $$\E_x[\tau_{A_M}] = \sum_{k = 1}^\infty \P(\tau_{A_M} > k) \leq Cn + \sum_{k=1}^\infty Cn \P(\tau_{A_M} > Cnk). $$
    Given $k \geq 1$, suppose $\tau_{A_M} > Cnk$. Either $\tau_{A_M}^M > k$ or $\tau_{A_M}^M < k$; in the latter case, this means that the first $Cnk$ steps of $(X_t)_{t \geq 0}$ do not hit $A_M$ but the first $k$ steps of the induced chain do, implying that there are $Cnk$ steps of $(X_t)_t$ that hit $\Omega_M$ fewer than $k$ times. Hence, $\tau_k > Cnk$. We then have
    \begin{align*} \sum_{k=1}^\infty \P(\tau_{A_M} > Cnk) &\leq Cn + \sum_{k=1}^\infty (\P(\tau_{A_M}^M > k) + \P(\tau_k > Cnk)) \\
    &= Cn + \E[\tau_{A_M}^M] + \sum_{k=1}^\infty \P(\tau_k > Cnk)\,.
    \end{align*}

    To bound $\E[\tau_{A_M}^M]$, observe that our upper bound on $t_M$ and \Cref{thm:hitting-time} imply that when we look at the induced chain, the expected number of $P_M$ transitions until the trajectory hits $A_M$ is $O(n^4 \log n)$.
    
    For the last term, let $\tau_{\Omega_M}$ be the hitting time of $\Omega_M$ in $(X_t)_t$. Note that $\min_{\sigma \in \Omega} \P(X_{t+1} \in \Omega_M \,|\, X_{t} = \sigma) \geq 1/n$ since $\sigma$ must have at least one MM face. Thus, $\E [\tau_{\Omega_M}] \leq n$. Conditioned on the $\sigma$-algebra generated by $\tau_{j-1}$, the random variable $Y_j = \tau_j - \tau_{j-1}$ is stochastically dominated by a geometric random variable of mean $n$.  If we let $\{\xi_j\}_{j\geq1}$ be independent and identically distributed geometric random variables of mean $n$, then we see as well that the sum $\tau_j = \sum_{i = 1}^j Y_j$ is stochastically dominated by $\sum_{i = 1}^j \xi_i$.   Thus,
    $$\P(\tau_k > Cnk) = \P\left(\sum_{j = 1}^k Y_j \geq C nk\right) \leq \P\left(\sum_{j = 1}^k \xi_k \geq C nk\right) \leq 2 \exp[-ck] $$
     for some constant $c > 0$, where the last inequality is by a general form of Bernstein's inequality for sub-exponential random variables (see, e.g., \cite{vershynin}).

    This shows $$\sum_{k = 1}^\infty \P(\tau_k > Cnk) = O(n)\,.$$

    Putting everything together, we then have
    $$\E_x[\tau_{A_M}] \leq O(n^5 \log n),$$
    and applying \Cref{thm:hitting-time} gives the claim.
\end{proof}

\section{Global Flat-Foldability}\label{sec:global}

\begin{figure}[h]
    \centering
    \includegraphics[scale=.45]{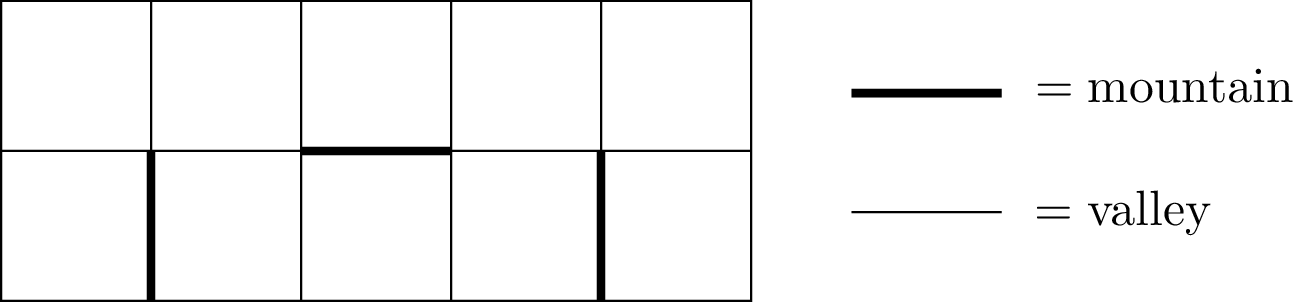}
    \caption{The non-flat-foldable MV assignment $\sigsp$ on the crease pattern $S_{2,5}$.\\}
    \label{fig:2x5}
\end{figure}

As previously mentioned, it is NP-hard to determine if a general origami crease pattern with a locally valid MV assignment is globally flat-foldable. This remains true even if we restrict ourselves to crease patterns whose vertices lie on a square lattice and all creases are only at angles of multiples of $45^\circ$ \cite{boxpleat}. Surprisingly, the same complexity problem for square grid crease patterns $S_{m,n}$ is unknown \cite{GFALOP,Origametry}. The smallest known example of a non-flat-foldable MV assignment on a square grid crease pattern was discovered by Justin in \cite{Justin} and is shown in Figure~\ref{fig:2x5}.

\begin{remark}
    The crease pattern in Figure~\ref{fig:2x5} is locally flat-foldable but not globally flat-foldable
\end{remark}

    A formal proof can be found in \cite[Example 7.12]{Origametry}. Informally, this can be seen by imagining that we fold the center horizontal mountain crease and the vertical valley creases to the immediate left and right of it. Folding these central creases forces the left- and right-most vertices of the crease pattern to be folded in the same direction, meaning that one would need to nest inside the other. Experimenting with an actual piece of paper shows, however, that such nesting is impossible to achieve.

\begin{figure}[h]
    \centering
    \begin{tikzpicture}[scale=0.7]
    \draw [step=1.0,black, very thick] (0,0) grid (5,4);  
    \draw [step=1.0,blue, very thick] (0,0) grid (4,3);
    \draw [step=1.0,red, very thick] (1,1) grid (3,2);
    \end{tikzpicture}
    \caption{An example of a neighborhood. Outlined in red is a $1 \times 2$ subset $T$; adding the creases highlighted in blue gives $N(T)$. Moreover, $N(N(T))$ is the entire crease pattern $S_{4,5}$.}
    \label{fig:nbhd}
\end{figure}
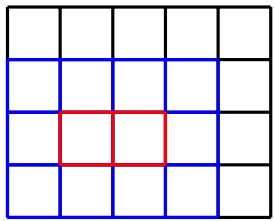

In this section, we aim to prove \Cref{thm:globalff-square}. Our first step is to show that we may extend a partial assignment of creases on a rectangular subset and the complement of its neighborhood in a way that preserves local flat-foldability.
Given an $a \times b$ subset $T$ of $S_{m,n}$, define the \textit{(closed) neighborhood} of $T$ as the crease pattern subset consisting of $T$ along with (when they exist) an additional column to the left and right of $T$ and an additional row above and below $T$. We notate this $N(T)$; see \Cref{fig:nbhd}.

\begin{lemma}\label{lem:extendingLFF}
    Let $a \leq m-2$ and $b \leq n-2$.
    Given locally flat-foldable configurations $\sig$ on $S_{m,n}$ and $\tau$ on $S_{a,b}$, there exists a configuration $\sig'$ on $S_{m,n}$ such that for an arbitrary $a \times b$ subset $T$ of $S_{m,n}$, we have
    \begin{enumerate}[(a)]
        \item $\sig'|_T \equiv \tau$,
        \item $\sig'|_{S_{m,n} \setminus N(T)} \equiv \sig_{S_{m,n} \setminus N(T)}$, and 
        \item $\sig'$ is locally flat-foldable.
    \end{enumerate}
\end{lemma}

\begin{proof}
    Define $\sig'$ by assigning the creases inside $T$ according to the configuration $\tau$ and assigning the creases outside $N(T)$ according to $\sigma$. Then $\sig'$ satisfies (a) and (b). Note that the only creases with no assignment are those in $N(T) \setminus T$. We will show that there is a valid MV assignment for these creases that maintains local flat-foldability at every vertex. 

   Given a crease pattern $S$, let $\partial_v(S)$ be the set of vertices on the boundary of $S$, meaning those in $S$ but incident to at least one edge not in $S$. Let $\partial_e(S)$ be the set of edges incident to exactly one vertex of $S$. Observe that for each $e \in \partial_e(T)$, $e$ must have one endpoint in $\partial_v(T)$ and the other endpoint in $\partial_v(N(T))$.

     Let $\partial_v(T) = \{v_1, v_2, \dots, v_{2(a-1)+2(b-1)}, u_1, u_2, u_3, u_4\}$ where $u_1, u_2, u_3, u_4$ are the vertices at the corners of $T$. For all $i$ such that $v_i \notin \partial_v(S)$, let $e_i = \{v_i, x_i\} \in \partial_e(T)$ where $x_i \in \partial_v(N(T))$.

    Observe that for every $i$, all creases incident to $v_i$ other than $e_i$ (if it exists) have an assignment induced by $\tau$. By Maekawa's Theorem (\ref{thm:maekawa}), then, there is a unique MV assignment to $e_i$ that maintains local flat-foldability. 

    The vertices $\{u_i : 1 \leq i \leq 4\} \cup \{x_i : 1 \leq i \leq 2(a+b-2), v_i \notin \partial_v(S)\}$ now each have two incident creases with an assignment and two incident creases without. It is straightforward to check that we may arbitrarily assign a subset of these creases (for example, one crease at each corner $u_i$ along with $x_1x_2$), and the remaining creases will have a unique assignment satisfying Maekawa's Theorem. 
\end{proof}

Looking towards proving Theorem \ref{thm:globalff-square}, we first prove an exponential upper bound on the probability of global flat-foldability:

\begin{lemma}\label{lem:globalff}
    Let $\sigma$ be a uniformly random locally flat-foldable configuration on $S_{m,n}$ where $m \geq 3, n \geq 6$. The probability that $\sigma$ is globally flat-foldable is at most $\exp[-cmn]$ for an absolute constant $c > 0$.
\end{lemma}

\begin{proof}
    Observe that if there exists $S_{a,b} \subset S_{m,n}$ such that $\sigma|_{S_{a,b}}$ is not globally flat-foldable, then $\sigma$ is also not globally flat-foldable. Consider the configuration $\sigsp$ given by Figure~\ref{fig:2x5}. 

    We may sample $\sig$ via the following procedure: partition $S_{m,n}$ into as many disjoint copies of $S_{3,6}$ as possible; label these $T_1, \dots, T_k$ where $k = \lfloor \frac{n}{6}\rfloor \cdot \lfloor \frac{m}{3}\rfloor$. Begin by sampling the creases of $T_1$ uniformly at random. Having sampled the faces of $T_1, \dots, T_i$, sample the unassigned creases of $T_{i+1}$ uniformly at random, conditioning on any neighboring creases that have already been assigned so that the local flat-foldability conditions are satisfied at each vertex. Note that by \Cref{lem:DLR}, this is equivalent to sampling $\sigma$ uniformly at random.
    
    Let $X_j$ be the indicator of the event that $\sigma$ restricted to the top-left $2 \times 5$ subset of $T_j$ is equivalent to $\sigsp$. Observe that $\P(X_j = 1) > 0$ due to Lemma~\ref{lem:extendingLFF}: given that $T_j$ contains $\sigsp$ in its upper-left corner, there exists a locally valid MV assignment to the remaining creases of $T_j$ that is compatible with arbitrary boundary conditions on $T_j$.

    Moreover, $\P(X_j = 1) = 2^{-10}$ so 
    $$\P(X_j \neq 1 \text{ for all } j) = (1-2^{-10})^k$$
    Thus, the probability that $\sig$ is globally flat-foldable is at most $(1-2^{-10})^k \leq e^{-cmn}$ where $c > 0$ is an absolute constant. 
\end{proof}

A submultiplicativity argument will then prove Theorem \ref{thm:globalff-square}.

\begin{proof}[Proof of Theorem \ref{thm:globalff-square}]
    Let $a_{m,n}$ be the number of globally flat-foldable configurations on $S_{m,n}$. Observe that $\{a_{m,n}\}_{m,n}$ is a submultiplicative sequence, meaning
    $$a_{m+m',n} \leq a_{m,n}a_{m',n}$$
    for all $m,m',n \in \mathbb N$. To see this, consider a globally flat-foldable configuration on $S_{m+m',n}$. The configurations induced on the leftmost $S_{m,n}$ and rightmost $S_{m',n}$ subgrids must be globally flat-foldable, and the inequality follows.

    By a multivariate analogue of Fekete's Lemma \cite{capobianco2008multidimensional}, we thus have that $\lim_{m,n \to \infty} \frac{\log a_{m,n}}{mn}$ exists and is equal to $\inf_{m,n \geq 1} \frac{\log a_{m,n}}{mn}$.  Letting $p_{m,n}$ denote the probability that a uniformly random locally flat-foldable configuration of $S_{m,n}$ is globally flat-foldable, this implies that $\lim_{m,n\to\infty} \frac{-\log p_{m,n}}{mn} =: c$ exists.
    By \Cref{lem:globalff} we in fact have that $c > 0$ thus completing the proof of the theorem.
\end{proof}

\begin{remark}
    Our proof for the square grid only required the existence of a single locally but not globally flat-foldable crease pattern and the ability to extend a partial locally flat-foldable crease pattern to a complete one, as in \Cref{lem:extendingLFF}. Thus, \Cref{thm:globalff-square} also holds for other lattices, such as the Miura-ori. 
\end{remark}

\section{Conclusion}\label{sec:conclusion}

We end this paper with comments on other classes of flat-foldable crease patterns, providing one example where we conjecture that face flips do give reasonable mixing times but proofs remain elusive,  and another where face flips do not give a reasonable way to randomly sample MV configurations.

\subsection{Triangle lattice}

\begin{figure}[h]
    \centering
    \includegraphics[scale=.6]{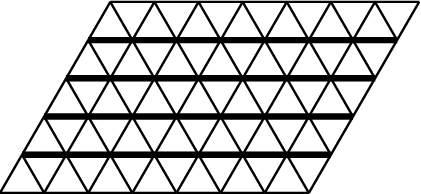}
    \caption{A triangular lattice crease pattern $T_{5,7}$ with a valid MV assignment.\\}
    \label{fig:triangular}
\end{figure}

Let the \textit{triangle lattice} crease pattern $T_{m,n}$ denote a parallelogram-shaped piece of paper with a crease pattern made entirely of tiled equilateral triangles, with $m$ and $n$ triangles along the two adjacent sides of the parallelogram. An example of a locally valid configuration for $T_{5,7}$ is shown in Figure~\ref{fig:triangular}. Unlike the square grid, the local and global flat-foldability of triangle lattice crease patterns are fairly unexplored.  
In \cite{faceflips} it is shown that the  origami flip graph of $T_{m,n}$ is connected with diameter $O(mn)$ and that
the problem of finding a shortest path between two vertices in $OFG(T_{m,n})$ is NP-complete. But these are the only known results for this class of crease patterns.  We note that considerable statistical physics work has been done on enumerating locally flat-folded crease patterns where one only folds a subset of edges of $T_{m,n}$ and  does not distinguish between mountain and valley \cite{DiF}.  This is a different notion of foldability than considered here and algorithmic questions  are open and interesting in that setting.  

\begin{figure}
    \centering
    \includegraphics[scale=.15]{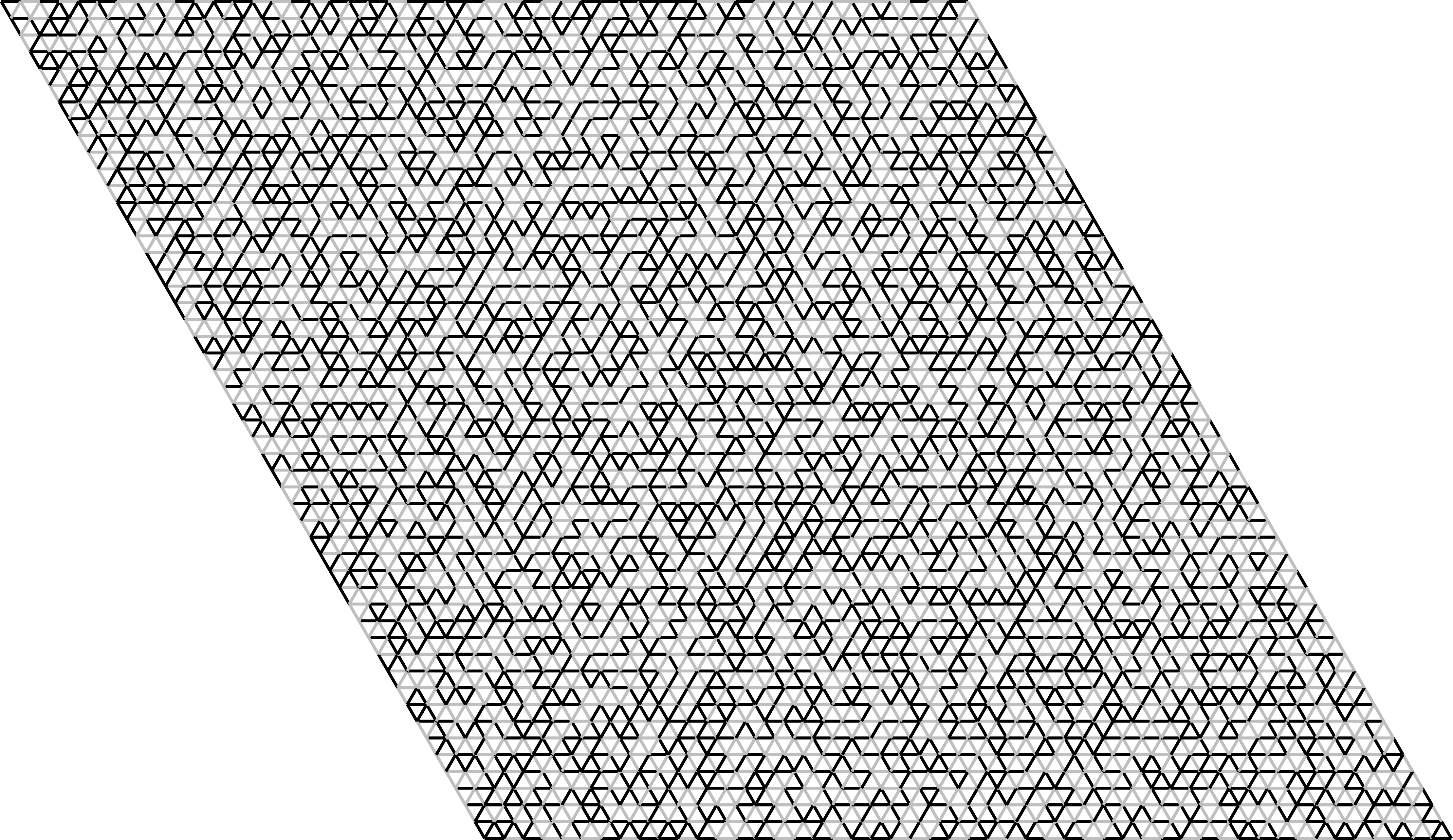}
    \caption{The results from a face-flip simulation, flipping the faces at random in $T_{50,50}$ to achieve a seemingly well-mixed result after $50^4$ iterations after starting from a configuration similar to \Cref{fig:triangular}. 
    }
    \label{fig:trimix}
\end{figure}

However, based on repeated simulations of the face-flip chain on relatively large examples---such as in \Cref{fig:trimix}---we conjecture that fast, and in fact optimal, mixing should occur.

\begin{conj}
    The mixing time of the face-flip Markov chain on $T_{m,n}$ is $O(mn \log(mn))$.
\end{conj}

\subsection{Complexity and the kite pattern}

The face-flip Markov chain is not the only avenue towards sampling algorithms. For example, there exist crease patterns for which no faces are flippable, e.g. the kite crease pattern, in which case a different approach has to be considered.

\begin{figure}
    \centering
\includegraphics[scale=.7]{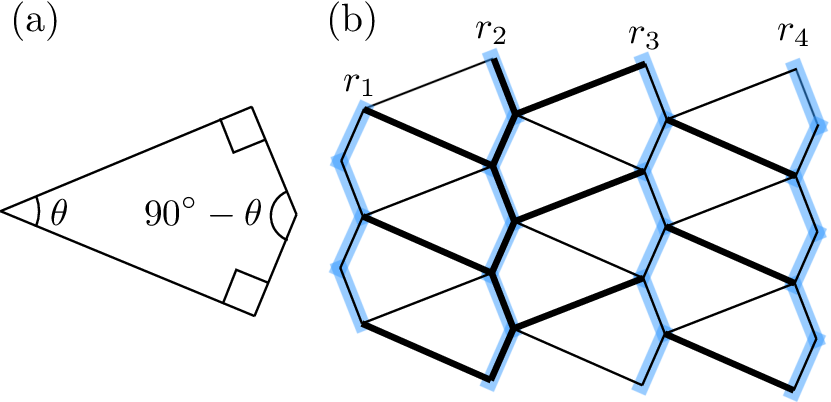}
    \caption{(a) A kite. (b) A tiling of kites with a MV assignment shown (mountains are bold lines, valleys are non-bold) and sets of same-configuration $r_i$ highlighted.}
    \label{fig:kite}
\end{figure}

As described in \cite{Assis}, and referred to as the \textit{Huffman grid} in \cite{faceflips}, the \textit{kite crease pattern} is made by tiling congruent quadrilaterals that have two opposite interior angles measuring $90^\circ$, with the other pair of angles being $\theta$ and $90^\circ-\theta$; see Figure~\ref{fig:kite}. As proven in \cite{faceflips}, any face flip of a locally valid MV assignment of this crease pattern results in a violation of the Big-Little-Big Theorem (\ref{thm:blb}). Therefore, none of the faces in a kite crease pattern are flippable and its origami flip graph has no edges.

However, there are many locally valid MV assignments of a kite crease pattern.  One is shown in Figure~\ref{fig:kite}(b), and as proven in \cite{faceflips}, each of the highlighted sets $r_i$ in the figure must be all-mountain or all-valley creases in order for the pattern to be valid. Thus, to get a different locally valid MV assignment, we could ``flip" any of these sets $r_i$ from all-mountain to all-valley, or vice-versa. Similarly, the sets of creases between consecutive sets $r_i$ and $r_{i+1}$ must alternate Ms and Vs, and any of these could be flipped, say from MVMVMV to VMVMVM.  Any of these kind of ``flips" of whole sets of creases can be made independently of the others and still result in a valid MV assignment. Therefore, if a kite crease pattern can be partitioned into $K$ of these flippable sets, then it will have exactly $2^K$ valid MV assignments. 
Thus we believe that similar arguments to Section~\ref{sec:sqtwist} will give results on the mixing time for kite crease patterns but without using the face-flip approach.

For the crease patterns we have studied so far, all of our results and conjectures have been in the positive direction---that there does indeed exist a computationally efficient way to sample locally flat-foldable assignments. But are there results in the negative direction for other equally ``natural'' crease patterns?

\begin{question}
    Are there natural families of crease patterns on which it is computationally hard to sample (either perfectly or approximately) from the uniform distribution on locally flat-foldable MV assignments?
\end{question}

A weaker question would be to only consider Markov chain sampling:

\begin{question}
    Are there natural families of crease patterns on which the face-flip Markov chain mixes slowly?
\end{question}

Now suppose that we are given a locally flat-foldable MV assignment $\sig_S$ to a subset $S$ of creases in a crease pattern $C$, and we ask if $\sig_S$ can be extended to a locally flat-foldable MV assignment for all of $C$. When $C$ is an $m\times n$ Miura-ori crease pattern, this question is answered in \cite{Ball} via an algorithm to decide the question in $O(mn\sqrt{mn})$ time.
In general, however, this problem remains open.

\begin{question}\label{question36}
    Given a locally flat-foldable crease pattern $C$ and a partial MV assignment, can we efficiently determine an MV assignment of the remaining creases that makes the entire assignment locally flat-foldable?
\end{question}

Lastly, the complexity of global flat-foldability remains a mystery. While the decision problem is NP-hard in general, the complexity for seemingly simple crease patterns remains unknown.

\begin{question}
    Is it computationally hard to sample a uniformly random globally flat-foldable MV assignment on the square grid?
\end{question}

\section*{Acknowledgments}
Thomas Hull is supported in part by NSF grants DMS-2428771 and DMS-2347000. Marcus Michelen is supported in part by NSF CAREER grant DMS-2336788 as well as DMS-2246624. Corrine Yap is supported in part by NIH grant R01GM126554 and in part by NSF grant DMS-1928930 while in residence at the Simons--Laufer Mathematical Sciences Institute during the Spring 2025 semester.

\bibliographystyle{plain}
\bibliography{references}

\begin{thebibliography}{10}

\bibitem{boxpleat}
Hugo~A. Akitaya, Kenneth~C. Cheung, Erik~D. Demaine, Takashi Horiyama,
  Thomas~C. Hull, Jason~S. Ku, Tomohiro Tachi, and Ryuhei Uehara.
\newblock Box pleating is hard.
\newblock In Jin Akiyama, Hiro Ito, Toshinori Sakai, and Yushi Uno, editors,
  {\em Discrete and Computational Geometry and Graphs}, pages 167--179, Cham,
  2016. Springer International Publishing.

\bibitem{faceflips}
Hugo~A. Akitaya, Vida Dujmovi{\'c}, David Eppstein, Thomas~C. Hull, Kshitij
  Jain, and Anna Lubiw.
\newblock Face flips in origami tessellations.
\newblock {\em Journal of Computational Geometry}, 7(1), 2016.

\bibitem{anari2021spectral}
Nima Anari, Kuikui Liu, and Shayan~Oveis Gharan.
\newblock Spectral independence in high-dimensional expanders and applications
  to the hardcore model.
\newblock {\em SIAM Journal on Computing}, (0):FOCS20--1, 2021.

\bibitem{Assis}
Michael Assis.
\newblock Exactly solvable flat-foldable quadrilateral origami tilings.
\newblock {\em Physical Review E}, 98(3):032112, Sep 2018.
\newblock arXiv: 1705.04710.

\bibitem{Ball}
Brad Ballinger, Mirela Damian, David Eppstein, Robin~Y. Flatland, Jessica
  Ginepro, and Thomas Hull.
\newblock {Minimum forcing sets for Miura folding patterns}.
\newblock In {\em Proc. 26th ACM-SIAM Symp. Discrete Algorithms, San Diego,
  California, 2015}, pages 136{--}147, 2015.

\bibitem{capobianco2008multidimensional}
Silvio Capobianco.
\newblock Multidimensional cellular automata and generalization of {F}ekete's
  lemma.
\newblock {\em Discrete Mathematics \& Theoretical Computer Science},
  10(Automata, Logic and Semantics), 2008.

\bibitem{chen2021optimal}
Zongchen Chen, Kuikui Liu, and Eric Vigoda.
\newblock Optimal mixing of {G}lauber dynamics: Entropy factorization via
  high-dimensional expansion.
\newblock In {\em Proceedings of the 53rd Annual ACM SIGACT Symposium on Theory
  of Computing}, pages 1537--1550, 2021.

\bibitem{GFALOP}
Erik~D. Demaine and Joseph O'Rourke.
\newblock {\em Geometric Folding Algorithms: Linkages, Origami, Polyhedra}.
\newblock Cambridge University Press, Cambridge, 2007.

\bibitem{DiF}
Philippe~Di Francesco.
\newblock Folding and coloring problems in mathematics and physics.
\newblock {\em Bulletin of the American Mathematical Society}, 37:251--307,
  2000.

\bibitem{Ginepro}
Jessica Ginepro and Thomas~C. Hull.
\newblock Counting {M}iura-ori foldings.
\newblock {\em Journal of Integer Sequences}, 17(10):Article 14.10.8, 2014.

\bibitem{Gjerde}
Eric Gjerde.
\newblock {\em Origami Tessellations: Awe-Inspiring Geometric Designs}.
\newblock {A K Peters}, Wellesley, MA, 2008.

\bibitem{Goldberg}
Leslie~Ann Goldberg, Russell Martin, and Mike Paterson.
\newblock Random sampling of 3-colorings in $\mathbb{Z}^2$.
\newblock {\em Random Structures \& Algorithms}, 24(3):279--302, 2004.

\bibitem{grimmett2018probability}
Geoffrey Grimmett.
\newblock {\em Probability on graphs: random processes on graphs and lattices},
  volume~8.
\newblock Cambridge University Press, 2018.

\bibitem{Active}
Edwin A.~Peraza Hernandez, Darren~J. Hartl, and Dimitris~C. Lagoudas.
\newblock {\em Active Origami: Modeling, Design, and Applications}.
\newblock Springer, 2019.

\bibitem{Tom1}
Thomas~C. Hull.
\newblock On the mathematics of flat origamis.
\newblock {\em Congressus Numerantium}, 100:215--224, 1994.

\bibitem{Tom2}
Thomas~C. Hull.
\newblock Counting mountain-valley assignments for flat folds.
\newblock {\em Ars Combinatoria}, 67:175--188, 2003.

\bibitem{Origametry}
Thomas~C. Hull.
\newblock {\em Origametry: Mathematical Methods in Paper Folding}.
\newblock Cambridge University Press, 2020.

\bibitem{JGAA-605}
{Thomas}~C. {Hull}, {Manuel} {Morales}, {Sarah} {Nash}, and {Natalya}
  {Ter-Saakov}.
\newblock Maximal origami flip graphs of flat-foldable vertices: properties and
  algorithms.
\newblock {\em Journal of Graph Algorithms and Applications}, 26(4):503--517,
  2022.
\newblock \url{https://jgaa.info/getPaper?id=605}.

\bibitem{Ida}
Tetsuo Ida.
\newblock {\em Introduction to Computational Origami}.
\newblock Springer, Cham, Switzerland, 2020.

\bibitem{Justin}
Jacques Justin.
\newblock Toward a mathematical theory of origami.
\newblock In Koryo Miura, editor, {\em Origami Science \& Art: Proceedings of
  the Second International Meeting of Origami Science and Scientific Origami},
  pages 15--29, Otsu, Japan, 1997. Seian University of Art and Design.

\bibitem{koehler}
J.~Koehler.
\newblock Folding a strip of stamps.
\newblock {\em Journal of Combinatorial Theory}, 5:135--152, 1968.

\bibitem{TessLang}
Robert~J. Lang.
\newblock {Tessellatica 11.1}.
\newblock Software, 2018.

\bibitem{levin2017markov}
David~A Levin and Yuval Peres.
\newblock {\em Markov chains and mixing times}, volume 107.
\newblock American Mathematical Soc., 2017.

\bibitem{Liu}
Bin Liu, Jesse~L. Silverberg, Arthur~A. Evans, Christian~D. Santangelo,
  Robert~J. Lang, Thomas~C. Hull, and Itai Cohen.
\newblock Topological kinematics of origami metamaterials.
\newblock {\em Nature Physics}, 14(8):811--815, 2018.

\bibitem{lunnon1}
W.~F. Lunnon.
\newblock A map-folding problem.
\newblock {\em Mathematics of Computation}, 22(101):193--199, 1968.

\bibitem{Maleczek2018}
Rupert Maleczek, Gabriel Stern, Clemens Preisinger, Moritz Heimrath, Oliver~D
  Krieg, and Astrid Metzler.
\newblock Curved folded wooden assemblies.
\newblock {\em Proceedings of IASS Annual Symposia}, 2018(20):1--8, 2018.

\bibitem{metropolis1953equation}
Nicholas Metropolis, Arianna~W Rosenbluth, Marshall~N Rosenbluth, Augusta~H
  Teller, and Edward Teller.
\newblock Equation of state calculations by fast computing machines.
\newblock {\em The Journal of Chemical Physics}, 21(6):1087--1092, 1953.

\bibitem{YOrigami}
Davic~C. Morgan, Denise~M. Halverson, Spencer~P. Magleby, Terri~C. Bateman, and
  Larry~L. Howell.
\newblock {\em Y Origami?: Explorations in Folding}.
\newblock American Mathematical Society, 2017.

\bibitem{nakajima2}
Chihiro Nakajima.
\newblock An efficient enumeration of flat-foldings : Study on random single
  vertex origami, 2024.

\bibitem{nakajima1}
Chihiro Nakajima.
\newblock A spin model for global flat-foldability of random origami, 2024.

\bibitem{peres2015mixing}
Yuval Peres and Perla Sousi.
\newblock Mixing times are hitting times of large sets.
\newblock {\em Journal of Theoretical Probability}, 28(2):488--519, 2015.

\bibitem{Robertson}
S.~A. Robertson.
\newblock Isometric folding of {R}iemannian manifolds.
\newblock {\em Proceedings of the Royal Society of Edinburgh},
  79(3--4):275--284, 1977--1978.

\bibitem{Silverberg1}
Jesse~L. Silverberg, Arthur~A. Evans, Lauren McLeod, Ryan~C. Hayward, Thomas
  Hull, Christian~D. Santangelo, and Itai Cohen.
\newblock Using origami design principles to fold reprogrammable mechanical
  metamaterials.
\newblock {\em Science}, 345(6197):647--650, 2014.

\bibitem{Silverberg}
Jesse~L. Silverberg, Jun-Hee Na, Arthur~A. Evans, Bin Liu, Thomas~C. Hull,
  Christian~D. Santangelo, Robert~J. Lang, Ryan~C. Hayward, and Itai Cohen.
\newblock {Origami structures with a critical transition to bistability arising
  from hidden degrees of freedom}.
\newblock {\em Nature Materials}, 14(4):389{--}393, 2015.

\bibitem{Uehara1}
Ryuhei Uehara.
\newblock {\em Introduction to Computational Origami}.
\newblock Springer Singapore, Singapore, 2020.

\bibitem{vershynin}
Roman Vershynin.
\newblock {\em High-dimensional probability: An introduction with applications
  in data science}, volume~47.
\newblock Cambridge university press, 2018.

\bibitem{wilson2004mixing}
David~Bruce Wilson.
\newblock Mixing times of lozenge tiling and card shuffling {M}arkov chains.
\newblock {\em The Annals of Applied Probability}, 14(1):274--325, 2004.

\bibitem{YouSci14}
Zhong You.
\newblock Folding structures out of flat materials.
\newblock {\em Science}, 345(6197):623--624, 2014.

\end{thebibliography}

\end{document}